\documentclass[12pt]{article}

\font\tenmsy=msbm10 scaled 1200
\font\sevenmsy=msbm7 scaled 1200
\font\fivemsy=msbm5 scaled 1200
\newfam\msyfam
\textfont\msyfam=\tenmsy
\scriptfont\msyfam=\sevenmsy
\scriptscriptfont\msyfam=\fivemsy
\newcommand{\Bbb}[1]{\fam\msyfam\relax#1}
\newcommand{\R}{{\Bbb R}}

\newcommand{\Z}{{\Bbb Z}}
\newcommand{\N}{{\Bbb N}}

\def\p{\psi}

\def\ep{ \epsilon }
\def\lp{\lambda_{\p}}
\def\D{{\rm D} }
\def\L{{\rm L} }

\newtheorem{theorem}{Theorem}
\newtheorem{corollary}{Corollary}

\newtheorem{lemma}{Lemma}

\newtheorem{thschmidt}{Covering Lemma }

\newtheorem{thbs}{Theorem KLW }

\newcommand{\B}{{\cal B}}

\begin{document}

\title{Metric Diophantine approximation and \\  `absolutely friendly' measures  }

\author{Andrew Pollington \footnote{This work was supported in part by the
National Science Foundation.}
 \\ {\small
B.Y.U., U{\scriptsize TAH}} \and Sanju Velani\footnote{Royal
Society University Research Fellow.} \\ {\small Y{\scriptsize
ORK}} }

\date{\small{For Ayesha and Iona on ${\rm N^o \, 2}$} \\ ~ \\
}

\maketitle

\abstract{Let $W(\p)$ denote the set of $\p$-well approximable
points in $\R^d$ and let $K$ be a compact subset of $\R^d$ which
supports a measure $\mu$.  In this short note, we show  that if
$\mu$ is an `absolutely friendly'  measure and a certain
$\mu$--volume sum converges then $\mu (W(\p) \cap K) = 0$. The
result obtained is in some sense analogous to the convergence part
of Khintchines classical theorem in the theory of metric
Diophantine approximation. The class of  absolutely friendly
measures is a subclass  of the friendly measures introduced in
\cite{KLW} and includes  measures supported on self similar sets
satisfying the open set condition. We also obtain an upper bound
result for the Hausdorff dimension of $W(\p) \cap K $. }

\section{Introduction}

{\small \subsection{The problem and results \ \  }}

The classical result of Dirichlet in the theory of Diophantine
approximation states that for any point  $ {\bf x}= (x_1, \dots,
x_d) \in \R^d $, there exist infinitely many $({\bf p},q) \in \Z^d
\times \N$ such that $$ \max_{1\leq i \leq d } |x_i -p_i/q | \leq
q^{-(d+1)/d}
\ \ \ . $$ Given a  real, positive  decreasing function $\psi :
\R^+ \to \R^+$,   a point ${\bf x} \in \R^d$ is said to be
$\psi$--well approximable if the above  inequality remains valid
with the right hand side replaced with $\psi(q)$.    We will
denote by $W(\p)$ the set of all such points; that is $$ W (\p) :=
\{{\bf x}  \in \R^d :\max_{1\leq i \leq d } |x_i - p_i/q| \leq
\p(q) \; {\rm for\ infinitely\ many\ } \; ({\bf p},q) \in \Z^d
\times \N \} \ . $$

\noindent A straightforward application of the Borel-Cantelli
lemma from probability theory yields the following statement.

\begin{lemma} $$|W(\p)|_d = 0 \hspace{7mm} {\rm if } \hspace{7mm}
\sum_{r=1}^{\infty} (r \, \p(r))^d \ < \ \infty \ \ .
$$\label{aim}
\end{lemma}
Thus, if the above sum converges then almost every (with respect
to $d$-dimensional Lebesgue measure) point ${\bf x } \in \R^d$ is
not $\p$--well approximable.  For $\tau \geq 0 $, consider the
function $\p_{\tau} : r \to r^{- \tau}$ and write  $W(\tau)$ for
$W(\p_\tau)$. In view of Dirichlet's result, $W(\tau)= \R^d$ for
$\tau \! \leq \! \mbox{\footnotesize $(d\!+\!1)/d $ }$. However,
in view of the above lemma we have that  $|W(\tau)|_d = 0 $ for
$\tau \! > \! \mbox{\footnotesize $(d\!+\!1)/d $ }$.

\vspace{6mm}

Now, let $K$ be a compact subset of $\R^d$ which supports a
non-atomic, finite  measure $\mu$ and let $$ W_K(\p) \ := \ K \cap
W(\p) \ \ . $$

\noindent In short, the problem is to determine conditions on
$\mu$ and $\psi$ under which  $\mu(W_K(\p)) =0$; i.e. $\mu$-almost
every point ${\bf x } \in \R^d$ is not $\p$--well approximable.
Note that $\mu(W_K(\p)) = \mu(W(\p))$ since $\mu$ is supported on
$K$. For the motivation behind the problem we refer the reader to
\cite{KLW,WT,W}.

\medskip

 In \cite{KLW},
Kleinbock, Lindenstrauss and Weiss introduce the notion of a
`friendly' measure and  show that if $\mu$ is  friendly then
$\mu(W_K(\tau)) =0$ for $\tau \! > \! \mbox{\footnotesize
$(d\!+\!1)/d $ }   \footnote{They actually prove their result in
the multiplicative framework.}$. They also show that the class of
friendly measures include (i) volume measures on non-degenerate
manifolds and (ii) measures supported on self similar sets
satisfying the open set condition. In full generality,  the
definition of friendly is rather technical and will not be
reproduce here -- see \S2 of \cite{KLW}.

\medskip

Our aim is to obtain a statement more  in line with Lemma
\ref{aim} which also implies that $\mu(W_K(\tau)) =0$ for  $\tau
\! > \! \mbox{\footnotesize $(d\!+\!1)/d $ }$. To achieve this we
impose conditions on $\mu$ which are stronger than those of
friendly. Nevertheless, measures supported on self similar sets
satisfying the open set condition are still included -- see
\S\ref{MEG}. Unfortunately, volume measures on non-degenerate
manifolds and not included.

\vspace{6mm}

Let $B(x,r)$ be a ball in $\R^d$ with centre $x$ and radius $r$.
The measure $\mu$ is said to be {\em doubling} if there exist
strictly positive constants $D$ and $r_0$ such that $$\mu(B(x,2r))
\ \leq \ D \, \mu(B(x,r)) \hspace{9mm} \forall \  x \in K
\hspace{3mm} \forall \ r < r_0 \ .  $$ The following notion of
`absolutely decaying' is essentially taken from \cite{KLW}.  Let
$\L$ denote a generic $(d-1)$--dimensional hyperplane of $ \R^d$
and let $\L^{(\ep)}$ denote its $\ep$-neighborhood. We say that
$\mu$ is {\em absolutely $\alpha$-decaying} if there exist
strictly positive constants $ C, \alpha, r_0  $ such that for any
hyperplane $\L$ and any $\ep > 0$ $$ \mu\left(B(x,r) \cap
\L^{(\ep)} \right) \ \leq \ C \, \left(\frac{\ep}{r}
\right)^{\alpha} \mu(B(x,r)) \hspace{7mm} \forall \  x \in K
\hspace{3mm}  \forall \ r < r_0 \ .  $$

\noindent In the case $d=1$, the hyperplane $\L$ is simply a point
$a \in \R$ and $\L^{(\ep)}$ is the ball $B(a,\ep)$ centred at $a$
of radius $\ep$.

\medskip

\noindent{\em Remark.  \ } Let $B(a,r)$ be a ball in $\R^d$. A
straightforward  geometric argument shows that if $\mu$ is
absolutely $\alpha$-decaying, then for any $\ep <
\mbox{\footnotesize $1/4$} $ \begin{equation} \mu\left(B(a,\ep r)
\right) \ \leq \ C \, 2^\alpha \, \ep^{\alpha} \ \mu(B(a,r))
\hspace{7mm} \forall \ a \in \R^d \hspace{4mm} \forall \ r < r_0 \
.  \label{winq}\end{equation} This essentially corresponds to the
condition on $\mu$ imposed in \cite{WT,W}. Note that in the case
$d=1$, condition (\ref{winq}) is equivalent to absolutely
$\alpha$-decay.

\vspace{6mm}

 \noindent{\bf Definition \ } A measure $\mu $ is
said to be {\em absolutely $\alpha$--friendly} if it is doubling
and absolutely $\alpha$-decaying.

\vspace{6mm}

We prove the following analogue of Lemma \ref{aim}.

\begin{theorem}
Let $K$ be a compact subset of $\R^d$ equipped with an absolutely
$\alpha$--friendly measure $\mu$. Then $$\mu (W_K(\p)) = 0
\hspace{6mm} {\rm if \ } \hspace{6mm} \sum_{r=1}^{\infty}
r^{\alpha\frac{d+1}{d} -1 } \p(r)^{\alpha}   \ < \ \infty \ \ . $$
\label{main}
\end{theorem}

\noindent{\em Remark.  \ } In the case when $d=1$, it is possible
to remove the condition that $\mu$ is doubling from the definition
of absolutely $\alpha$--friendly; i.e. all that is required is
that $\mu$ is absolutely $\alpha$--decaying -- see \S\ref{d=1}.
With  this  in mind, the above theorem restricted to $d=1$  is
identical to that established in \cite{WT}.  Thus, Theorem
\ref{main} constitutes the natural higher dimensional analogue of
\cite{WT}. The above theorem should also be compared with Theorem
9 of \cite{W}.

\medskip


  Note that in the case that   $\p_{\tau} : r \to
r^{-\tau} $ and $\tau >  \mbox{\small $\frac{d+1}{d}$}$,

$$ \sum_{r=1}^{\infty} r^{\alpha\frac{d+1}{d} -1 }
\p_{\tau}(r)^{\alpha} \ := \ \sum_{r=1}^{\infty} r^{-1-\alpha(\tau
- \mbox{\tiny $\frac{d+1}{d}$})}
 \ < \ \infty $$ and so Theorem \ref{main} implies that
$\mu(W_K(\tau)) =0$ whenever $\mu$ is absolutely
$\alpha$--friendly. More to the point, consider the function $\p:
r \to r^{-\mbox{\tiny $\frac{d+1}{d}$} } (\log r)^{-\beta} $ where
$\beta > 1/\alpha $. Then $$ \sum_{r=1}^{\infty}
r^{\alpha\frac{d+1}{d} -1 } \p(r)^{\alpha} \ := \
\sum_{r=1}^{\infty} r^{-1}(\log r)^{-\alpha \beta} \ < \ \infty  \
, $$ and  Theorem \ref{main} implies that $\mu(W_K(\p)) =0$
whenever $\mu$ is absolutely $\alpha$--friendly.

\medskip

It will be evident from the proof that all that is actually
required in establishing the theorem is that the doubling and
absolutely $\alpha$-decaying inequalities  are satisfied at
$\mu$--almost every point in $ K$. Also the relevance of
hyperplanes  in the definition of absolutely $\alpha$-decaying
will become crystal clear from our proof of the theorem.
Essentially, on the real line $\R$ an interval $I_n$ of length
$1/4n^2$ can contain at most one rational $p/q$ with $n \leq q <
2n$. This follows from the trivial observation that if $n \leq q,
q' < 2n$ then  $|p/q-p'/q'| \geq 1/qq' > 1/ 4n^2$; i.e. the
distance between two such rationals is strictly greater than the
length of $I_n$. The higher dimension analogue of this is the
following. Let $B_n$ be a ball in $\R^d$ of radius $c/n^{(d+1)/d}$
where $c$ is a sufficiently small constant dependent only on $d$.
Then any rational points ${\bf p}/q$ lying within $B_n$ with
$n\leq q < 2n$ must lie on a single $(d\!-\!1)$--dimensional
hyperplane $\L$. This is the key observation on which the proofs
of Theorems \ref{main} and \ref{hdmain}  hinge.

\vspace{6mm}

We now turn our attention to determining an upper bound for $\dim
W_K(\p)$ -- the Hausdorff dimension of $W_K(\p)$. For $s \geq 0$,
let ${\cal H}^s$ denote the $s$--dimensional Hausdorff measure -
see \S\ref{hmd}.

\begin{theorem}
Let $K$ be a compact subset of $\R^d$ equipped with an absolutely
$\alpha$--friendly measure $\mu$. Furthermore, suppose  there
exist positive constants $a,b, \delta$ and $r_o$   such that
\begin{equation}
 a \, r ^{\delta}  \ \leq  \  \mu(B(x,r))  \ \leq  \   b \, r
^{\delta}  \hspace{7mm} \forall \  x \in K \hspace{3mm}  \forall \
r < r_0 \ .   \label{measure}
\end{equation}
 Then, for $ s \leq \delta$ $$ {\cal H}^s
(W_K(\p)) = 0 \hspace{6mm} {\rm if \ } \hspace{6mm}
\sum_{r=1}^{\infty} r^{\alpha\frac{d+1}{d} -1 } \p(r)^{\alpha+s -
\delta}   \ < \ \infty \ \ . $$

\label{hdmain}

\end{theorem}

\noindent{\em Remark 1: \ } Note that (\ref{measure})  imposed on
$\mu$ trivially implies that $\mu$ is doubling. Furthermore, if
$\delta \!> \!d\!-\!1$ then (\ref{measure}) together with a
straightforward geometric argument implies that $\mu$ is
absolutely $\alpha$--decaying with $\alpha :=\delta -(d-1) > 0 $.
Thus, if $\delta \!> \!d\!-\!1$ the hypothesis that $\mu$ is
absolutely $\alpha$--friendly is in fact  redundant from the
statement of Theorem \ref{hdmain}.

\medskip

\noindent{\em Remark 2: \ }If $K$ supports a measure $\mu$
satisfying (\ref{measure}), then   $\dim K = \delta$ and moreover
that $ 0< {\cal H}^\delta  (K) < \infty$ -- see \cite{falc} for
the details. Now, since $W_K(\p)$ is a subset of  $K $ we have
that $\dim W_K(\p) \leq \delta $ and so  ${\cal H}^s (W_K(\p)) = 0
$ for any $s > \delta$. Thus, the condition  $s \leq \delta$ in
the statement of the theorem can be assumed without any loss of
generality.

\medskip

Given a real, positive decreasing  function $\p$, the {\em lower
order} $\lp$ of $1/\p$ is defined by $$ \lp \ := \
\liminf_{r\to\infty}\; \frac{- \log \p(r)}{\log r} \ , $$ and
indicates the growth of the function $1/\psi$ `near' infinity.
Note that $\lp $  is non-negative since $\psi$ is a decreasing
function. A simple consequence of Theorem \ref{hdmain} is the
following statement.

\begin{corollary}
Let $K$ be a compact subset of $\R^d$ equipped with an absolutely
$\alpha$--friendly measure $\mu$ satisfying (\ref{measure}). Then,
for  $\lp \! \geq \! \mbox{\footnotesize $(d\!+\!1)/d $ }$ $$ \dim
W_K(\p) \ \leq \ \delta \, - \, \alpha \left(1 - \mbox{\small
$\frac{d+1}{\lp \,  d}$} \right)  \ \ . $$ \label{hdmaingencor}
\end{corollary}

As a special case we obtain the following statement.

\begin{corollary}
Let $K$ be a compact subset of $\R^d$ equipped with an absolutely
$\alpha$--friendly measure $\mu$ satisfying (\ref{measure}). Then,
for  $\tau \! \geq \! \mbox{\footnotesize $(d\!+\!1)/d $ }$ $$
\dim W_K(\tau) \ \leq \ \delta \, - \, \alpha \left(1 -
\mbox{\small $\frac{d+1}{\tau \,  d}$} \right)  \ \ . $$
\label{hdmaincor}
\end{corollary}

Note that for $\tau \! > \! \mbox{\footnotesize $(d\!+\!1)/d $ }$
we have that $\dim W_K(\tau) < \delta $. Since $\mu$ is comparable
to ${\cal H}^\delta $ restricted to $K$, it follows that
$\mu(W_K(\tau)) =0$.

\vspace{4mm}

\noindent{\em A general remark: \ }For $d \geq 2 $, it is highly
unlikely that either Theorem \ref{main} or Theorem \ref{hdmain}
are ever sharp.  For instance, take the case that $K := [0,1]^d $
and $\mu$ is $d$--dimensional Lebesgue measure.  It is easily
verified that $\mu$ is absolutely $\alpha$--friendly with $\alpha
=1 $.  Thus, Theorem \ref{main} implies that $|W_K(\p)|_d = 0$
whenever $$ \sum_{r=1}^{\infty} r^{\mbox{\small
$\frac{d+1}{d}$}-1} \, \p(r) \ < \ \infty \ \ . $$ So when $d=1$
this coincides with the Lemma \ref{aim}.  However, for $d \geq 2 $
the above statement is weaker than that of the lemma. In view of
Khintchines theorem one knows that the lemma is sharp; that is to
say that if the sum in the lemma diverges then not only is
$|W_K(\p)|_d > 0$ but it is of full measure. It is probable that
the theorems of this paper are sharp in the case $d=1$.

{\small  \subsection{The main example \ \ \label{MEG} }}

The following statement which combines  Theorems 2.2 and 8.1  of
\cite{KLW}, shows that a large class of fractal measures are
absolutely $\alpha$--friendly and satisfy (\ref{measure}).

\vspace{7mm}

\begin{thbs} Let $\{{\bf S}_1, \dots, {\bf S}_k \} $ be an
irreducible family of contracting self similarity  maps of $\R^d$
satisfying the open set condition  and let   $\mu$ be the
restriction of  ${\cal H}^{\delta}$ to its attractor $K$ where $
\delta := \dim K$.  Then $\mu$ is absolutely $\alpha$--friendly
and  satisfies (\ref{measure}).
\end{thbs}

\medskip

Thus for the natural measures associated with self similar sets
satisfying the open set condition, Theorems \ref{main} and
\ref{hdmain} are applicable. The simpliest  examples of such sets
include regular Cantor sets, the Sierpi$\acute{{\rm n}}$ski gasket
and the von Kock curve.  All the terminology except for
`irreducible' is pretty much standard -- see for example
\cite[Chp.9]{falc}. The notion of irreducible introduced in
\cite[\S2]{KLW}  avoids the natural obstruction that there is a
finite collection of proper affine subspaces of $\R^d$ which is
invariant under $\{{\bf S}_1, \dots, {\bf S}_k \} $.

\section{Hausdorff measures and dimension \label{hmd} }
 In this short section  we define Hausdorff measure and dimension
for completeness and in order to establish some notation. For
$\rho
> 0$, a countable collection $ \left\{B_{i} \right\} $ of
Euclidean balls in  $\R^{d}$ of radii $r_i \leq \rho $ for each
$i$ such that $X \subset \bigcup_{i} B_{i} $ is called a $ \rho
$-cover for $X$.  Let $s$ be a non-negative number and define $$
 {\cal H}^{s}_{ \rho } (X)
  \; = \; \inf \left\{ \sum_{i} r_i^s
\ :   \{ B_{i} \}  {\rm \  is\ a\  } \rho {\rm -cover\  of\ } X
\right\} \; , $$ where the infimum is taken over all possible $
\rho $-covers of $X$. The {\it s-dimensional Hausdorff measure}
${\cal H}^{s} (X)$ of $X$ is defined by $$ {\cal H}^{s} (X) =
\lim_{ \rho \rightarrow 0} {\cal H}^{s}_{ \rho } (X) = \sup_{ \rho
> 0} {\cal H}^{s}_{ \rho } (X)
$$ \noindent and the {\it Hausdorff dimension} dim $X$ of $X$ by
$$ \dim \, X = \inf \left\{ s : {\cal H}^{s} (X) =0 \right\} =
\sup \left\{ s : {\cal H}^{s} (X) = \infty \right\} \, . $$

\vskip 9pt

Further details and alternative definitions of Hausdorff measure
and dimension can be found in \cite{falc}.

\section{A covering lemma}

\vskip 4pt

The following rather simple covering result will be used at
various stages during  the proof of our theorems.

\medskip

\begin{thschmidt}
Let $(\Omega,d)$ be a  metric space and $\B$ be a finite
collection of balls with common radius $r>0$. Then there exists a
disjoint sub-collection $\{B_i\}$ such that $$\bigcup_{B\in\B} \,
B \ \subset \ \bigcup_i \  {3B}_i \ . $$
\end{thschmidt}

\medskip

\noindent{\em Proof : \ }Let $S$ denote the set of centres of the
balls in $\B$.  Choose $c_1 \in S$ and for $k \geq 1 $, $$c_{k+1}
\ \in \ S \ \backslash \  \bigcup_{i=1}^{k} B(c_i, 2r) $$ as long
as $S  \ \backslash \ \bigcup_{i=1}^{k} B(c_i, 2r) \neq \emptyset
$. Since $\#S$ is finite, there exists $k_1 \leq \#S $ such that
$$ S \ \subset \  \bigcup_{i=1}^{k_1} B(c_i, 2r) \ . $$ By
construction, any ball $B(c,r)$ in the original collection $\B$ is
contained in some ball $B(c_i,3r)$ and since $d(c_i,c_j) > 2r$
 the chosen balls $ B(c_i,r) $ are clearly  disjoint.

\hfill $ \spadesuit $

\section{Proof of Theorem \ref{main} \label{d2}  }

\noindent{{\bf Step 1}: {\em Preliminaries.} \ }  We are assuming
that  $\sum r^{\alpha\frac{d+1}{d} -1 } \p(r)^{\alpha} $
converges and since  $\p$ is monotonic, it follows that
\begin{equation}
\sum_{n=1}^{\infty} \left(2^{n\frac{d+1}{d} }
\p(2^n)\right)^{\alpha} < \infty   \  . \label{div2}
\end{equation}

\noindent Next notice, that without loss of generality we can
assume that
\begin{equation}
\p(2^n) \ < \ c \, 2^{- n\frac{d+1}{d} } \label{cpsi}
\end{equation} for any $c > 0 $
and $n$ sufficiently large.  This is easy to see. Suppose on the
contrary that  there exists a sequence $\{n_i\} $ such that
$\p(2^{n_i})  \geq  c \, 2^{- n_i \frac{d+1}{d} }  $. Then $$
\sum_{n=1}^{\infty}  \left(2^{n\frac{d+1}{d} }
\p(2^n)\right)^{\alpha} \ \geq \ \sum_{i=1}^{\infty}
\left(2^{n_i\frac{d+1}{d} } \p(2^{n_i})\right)^{\alpha} \ \geq \
c^{\alpha} \, \sum_{n=1}^{\infty}  1   \ = \ \infty  \ \ , $$ and
this contradicts (\ref{div2}).

\medskip

\noindent{{\bf Step 2}: {\em The balls $\D_n$.} \ } For $n \in
\N$, let $\D_n$ denote a generic ball with centre in $K$ and of
radius $$r_n := \mbox{\normalsize $ \frac{1}{6}$}
 \, \mbox{\normalsize $ \left(\frac{1}{\kappa \, d !
}\right)^{\frac{1}{d}} $} \  2^{- \frac{d+1}{d}(n+1) } \ \ . $$
Here $\kappa := \kappa(d) $ is the volume ($d$--dimensional
Lebesgue measure) of a ball of radius one in $\R^d$. In view of
the covering lemma and the fact that $K$ is compact, there exists
a finite, disjoint collection ${\cal D}_n $ of balls $\D_n$ with
centers in $K$ such that $$ \bigcup_{{\cal D}_n} 3 \D_n \ \supset
\ K  \  \ . $$

\noindent Note that since $\mu $ is doubling, we have that
\begin{equation}
\sum_{{\cal D}_n} \mu(3 \D_n)  \ \leq \ \sum_{{\cal D}_n} \mu(4
\D_n)  \ \leq \ D^2 \ \sum_{{\cal D}_n} \mu( \D_n) \ = \ D^2 \,
\mu\!\!\left(\bigcup_{{\cal D}_n}^{\circ} \D_n \right) \ \leq \
D^2 \, \mu(K) \ \ . \label{sumDn}
\end{equation}

\medskip

Next, consider a ball $3\D_n $ where $\D_n \in {\cal D}_n $.
Suppose there is a rational point ${\bf p}/q:=(p^{(1)}/q, \ldots,
p^{(d)}/q) $ such that
\begin{equation}
B({\bf p}/q, \sqrt d \, \psi(q))  \cap \ 3\D_n \neq \emptyset
\hspace{7mm} {\rm and \ } \hspace{7mm}  2^n \leq q < 2^{n+1} \ .
\label{badrat}
\end{equation}
By (\ref{cpsi}) and using the  fact that $\p$ is decreasing, it
follows that for $n$ sufficiently large ${\bf p}/q \in 6\D_n $.
Now assume that there are $d+1$ or more such rational points
satisfying (\ref{badrat}).  Take any $d+1$ such rationals; ${\bf
p}_{\atop{\mbox{\tiny $\!\!\! 0$}}}/q_{\atop{\mbox{\tiny $\!
0$}}}, {\bf p}_{\atop{\mbox{\tiny $\!\!\!
1$}}}/q_{\atop{\mbox{\tiny $\! 1$}}} ,\ldots , {\bf
p}_{\atop{\mbox{\tiny $\!\!\! d $}}}/q_{\atop{\mbox{\tiny $\! d
$}}} $.
 In view of the denominator
constraint, the rational points are necessarily distinct. Suppose
for the moment that they do not lie on a $(d\!-\!1)$--dimensional
hyperplane and form the $d$--dimensional simplex  $ \Delta$
sub-tended by them; i.e. an interval when $d=1$, a triangle when
$d=2$ and a tetrahedron when $d=3$. The volume ($d$--dimensional
Lebesgue measure) of the simplex $ \Delta$ times $d$ factorial  is
equal to the absolute value of the determinant $$ \det \ := \
\left|
\begin{array}{ccccc} 1 & p_{\atop{\mbox{\tiny $\!
0$}}}^{\mbox{\tiny (1)}}/q_{\atop{\mbox{\tiny $\! 0$}}} & \ldots &
p_{\atop{\mbox{\tiny $\! 0$}}}^{\mbox{\tiny
(d)}}/q_{\atop{\mbox{\tiny $\! 0$}}} \\  \\ 1 &
p_{\atop{\mbox{\tiny $\! 1$}}}^{\mbox{\tiny
(1)}}/q_{\atop{\mbox{\tiny $\! 1$}}} & \ldots &
p_{\atop{\mbox{\tiny $\! 1$}}}^{\mbox{\tiny
(d)}}/q_{\atop{\mbox{\tiny $\! 1$}}} \\
\\  \vdots \\  \\ 1 & p_{\atop{\mbox{\tiny $\!
d$}}}^{\mbox{\tiny (1)}}/q_{\atop{\mbox{\tiny $\! d$}}}  & \ldots
& p_{\atop{\mbox{\tiny $\! d$}}}^{\mbox{\tiny
(d)}}/q_{\atop{\mbox{\tiny $\! d$}}}
\end{array} \right|  \ . $$

\noindent Then, by (\ref{badrat}) $$ d !  \ \times  \  |\Delta|_d
\ \geq \ \frac{1}{q_{\atop{\mbox{\tiny $\! 0$}}}
q_{\atop{\mbox{\tiny $\! 1$}}} \ldots  q_{\atop{\mbox{\tiny $\!
d$}}}}
> 2^{-(d+1)(n+1)} \ . $$

\noindent Trivially, $$ |\,  6\D_n \,  |_d   \ =  \ \kappa \,  (6
\, r_n)^d \ := \ \mbox{\normalsize $ \frac{1}{d!} \ $}
2^{-(d+1)(n+1)}   \ \ . $$

\noindent Thus $ |\Delta |_d
> | 6\D_n |_d$ and this is impossible since
 $\Delta \subset 6\D_n$.  The upshot of this is that the $d$--dimensional simplex $\Delta $
can not exist and so if there are $d +1$ or more rational points
satisfying (\ref{badrat}) then they must lie on a
$(d\!-\!1)$--dimensional hyperplane  $\L := \L(\D_n)$ passing
through the ball $3\D_n$.  In the event that there are no more
than   $d$ rational points satisfying (\ref{badrat}), the
existence  of such a  hyperplane  is obvious -- of course it is
not  unique if the number of rational points is less than $d$.
Thus, associated with  each ball $\D_n \in  {\cal D}_n$ there is a
$(d\!-\!1)$--dimensional hyperplane  $\L := \L(\D_n)$ containing
all rational points satisfying (\ref{badrat}). Note that in the
case $d=1$, any hyperplane $\L$ is simply a point.

\medskip

\noindent{{\bf Step 3}: {\em The finale}. \ } For $n \in \N$, let
$$A_n \ := \ \bigcup_{2^n \leq q < 2^{n+1}} \bigcup_{{\bf p} \in
\Z^d} B \left({\bf p}/q, \sqrt d  \;  \p(q) \right)
 \ . $$
By definition, $W_K(\p) \subset \limsup_{n \to \infty} A_n \cap \,
K $. It follows via Step 2 and the fact that $\psi$ is decreasing,
that for $n$ sufficiently large
\begin{eqnarray*} \mu (A_n) \ := \ \mu (A_n \cap K ) & = &
\mu\left( A_n \cap \bigcup_{{\cal D}_n} 3 \D_n \right) \\ & \leq &
\sum_{{\cal D}_n} \mu\left(3 \D_n \cap \L^{(\ep)} \right)
\hspace{19mm} \ep := \sqrt d \, \p(2^n) \ \ \   \L := \L(\D_n)
\\ & \ll &  \left(2^{n\frac{d+1}{d} }
\p(2^n)\right)^{\alpha} \ \sum_{{\cal D}_n} \mu(3 \D_n)
\hspace{5mm}   \mu \mbox{ \ is absolutely $\alpha$--decaying} \\ &
\ll & \left(2^{n\frac{d+1}{d} } \p(2^n)\right)^{\alpha}  \ \mu(K)
\hspace{14mm}  {\rm by \ } (\ref{sumDn}).
\end{eqnarray*}

\noindent Hence, by (\ref{div2})

$$ \sum \mu(A_n \cap K ) \ = \  \sum \mu(A_n) \ \ll \  \sum
\left(2^{n\frac{d+1}{d} } \p(2^n)\right)^{\alpha} \ < \ \infty \ \
$$
\medskip

\noindent and the Borel-Cantelli lemma implies that $
\mu(\limsup_{n \to \infty}  A_n ) = 0 $. Thus,  $ \mu(W_K(\p))$ is
zero  as required.

\hfill $\spadesuit$

\subsection{The case when $d = 1$ revisited \label{d=1}}
Clearly the above proof contains the  case when $d=1$. However, it
is possible to give a more direct proof of a stronger statement
which does not assume that $\mu$ is doubling -- see the remark
straight after the statement of Theorem \ref{main}. Although  the
proof below is basically the same as that in \cite{WT}, we have
decided to include a sketch in order to bring out the true nature
of the `simplex/determinate' argument and the role of hyperplanes
when $d \geq 2 $ in the proof above. In the $d=1$ case, the
`simplex/determinate' argument reduces to the following. Consider
rationals $p/q$ with $2^n \leq q <2^{n+1}$. For any two such
rationals, notice that $$ \left|\frac{p}{q} - \frac{p'}{q'}
\right| \ \geq \ \frac{1}{qq'} \
> \ 2^{-2(n+1)} \ := \ 2 \, r_n \  . $$
Thus, any interval of length $2r_n$ can contain at most one
rational. In particular,  the intervals $B(p/q,r_n)$ are disjoint.

Now let $A_n$ be as in Step 3. By definition, $W_K(\p) =
\limsup_{n \to \infty} A_n \cap \, K $. Then, in view of
(\ref{cpsi}) and the fact that $\psi$ is decreasing  we have that
for $n$ sufficiently large
\begin{eqnarray*} \mu (A_n)  & \leq &
\!\!\!\! \sum_{2^n \leq q < 2^{n+1}} \sum_{p \in \Z} \mu(B(p/q,
\ep \, r_n))  \hspace{39mm} \ep \, r_n := \p(2^n)  \\ & & ~
\\ & \ll & \left(2^{2n}
\p(2^n)\right)^{\alpha} \  \sum_{2^n \leq q < 2^{n+1}} \sum_{p \in
\Z} \mu(B(p/q, r_n))  \hspace{16mm} {\rm by \ } (\ref{winq}) \\
 & \leq & \left(2^{2n}
\p(2^n)\right)^{\alpha} \ \mu(K) \hspace{52mm} {\rm by \ }
 \hfill {\rm \ disjointness.}
\end{eqnarray*}

The Borel-Cantelli lemma implies the desired statement. \hfill
$\spadesuit$

\section{Proof of Theorem \ref{hdmain} }

To a certain extent the proof of Theorem \ref{hdmain} is similar
to that of Theorem \ref{main}.

\medskip

\noindent{{\bf Step 1}: {\em Preliminaries.} \ } Without loss of
generality we can assume that $\p(r) \to 0 $ as $r \to \infty$.
Suppose that this was not the case. Then $W_K(\p) = K $ by
Dirichlet's theorem and so ${\cal H}^s (W_K(\p)) > 0 $ for any $s
\leq \delta$ - see Remark 2 straight after the statement of
Theorem \ref{hdmain}.

\medskip

Without loss of generality we can assume that $s > \delta - \alpha
$. If this where not the case then the sum in the statement of the
theorem cannot possibly converge.

\medskip

Since  $\p$ is monotonic, the convergence of the sum in the
statement of the theorem is equivalent to
\begin{equation}
\sum_{n=1}^{\infty} 2^{n \alpha \frac{d+1}{d} } \p(2^n)^{\alpha+s
- \delta} < \infty \ . \label{div2hd}
\end{equation}

\medskip

 Finally, notice that since $s > \delta - \alpha$,  we
can assume (\ref{cpsi})  without any loss of generality.
Otherwise, (\ref{div2hd}) would be contradicted.

\medskip

\noindent{{\bf Step 2}: {\em A good $\rho$--cover for $W_K(\p)$.}
\ }  For $n \in \N$, let ${\cal D}_n $ be the disjoint collection
of balls  $\D_n $ as defined in Step 2 of \S\ref{d2}. Since the
collection is disjoint and $\mu$ satisfies (\ref{measure}), we
have that for $n$ sufficiently large $$\#{\cal D}_n \times 2^{-
\mbox{\footnotesize $ \frac{d+1}{d}$} (n+1)\delta } \ \asymp \
\sum_{{\cal D}_n} \mu( \D_n) \ = \ \, \mu\!\!\left(\bigcup_{{\cal
D}_n}^{\circ} \D_n \right) \ \leq \ \mu(K) \ \ . $$ Thus, for $n$
sufficiently large
\begin{equation}
\#{\cal D}_n  \ \ll \ 2^{\mbox{\footnotesize $ \frac{d+1}{d}$}
(n+1)\delta }  \ \ . \label{nod}
\end{equation}

\medskip

Now put $\ep := \sqrt d \,  \p(2^n)$ and fix some ball $\D_n \in
{\cal D}_n $. Let  $ \L := \L(\D_n) $  be the   $(d \!-\!
1)$--dimensional hyperplane  associated with $\D_n$ -- see Step 2
of \S\ref{d2}. In view of the covering lemma, there exists a
finite disjoint collection ${\cal C }(\D_n)$ of balls $B_n(\p)$
with centers in $K$ and common radius $\p(2^n)$ such that
\begin{equation}
\bigcup_{{\cal C }(\D_n)} 3 B_n(\p) \ \supset \ 3 \, \D_n \cap
\L^{(\ep)}  \cap K \label{sanju}
\end{equation}
and
$$\bigcup_{{\cal C }(\D_n)}^{\circ}  B_n(\p) \ \subset \ 6 \, \D_n
\cap \L^{(2\ep)}
$$ The latter together with the fact that $\mu$ is absolutely
$\alpha$--decaying implies that
\begin{eqnarray*} \#{\cal C}(\D_n) \times \p(2^n)^{\delta } &
\asymp & \sum_{{\cal C}(\D_n)} \mu( B_n(\p)) \ = \ \,
\mu\!\!\left(\bigcup_{{\cal C}(\D_n)}^{\circ} B_n(\p) \right) \\ &
& \\ & \leq & \mu( 6 \, \D_n \cap \L^{(2\ep)}) \ \ll \ \left(
2^{\mbox{\footnotesize $ \frac{d+1}{d}$} (n+1)} \p(2^n)
\right)^{\alpha} \ 2^{- \mbox{\footnotesize $ \frac{d+1}{d}$}
(n+1)\delta }  \ .
\end{eqnarray*}
Thus, for $n$ sufficiently large
\begin{equation}
\#{\cal C}(\D_n) \ \ll \ \left( 2^{\mbox{\footnotesize $
\frac{d+1}{d}$} (n+1)} \p(2^n) \right)^{\alpha - \delta} \ \ .
\label{noc}
\end{equation}

\noindent Now with $A_n$ defined as in Step 3 of \S\ref{d2}, it
follows via  (\ref{sanju}) that
\begin{eqnarray*}A_n \cap K \  & = & \ \bigcup_{\D_n \in {\cal D}_n }
\!\!\!\! 3\D_n  \, \cap A_n  \cap K \  \ \subset  \ \
\bigcup_{\D_n \in {\cal D}_n } 3\D_n \cap \L^{(\ep)} \cap K  \\ &
&  \\ \ & \subset & \ \bigcup_{\D_n \in {\cal D}_n } \ \
\bigcup_{B_n(\p) \in {\cal C}(\D_n)} \!\! 3\, B_n(\p) \ \ .
\end{eqnarray*}

 In particular, for each $k \in \N$ the collection $$ \left\{
3\, B_n(\p) \; : \; B_n(\p) \in {\cal C}(\D_n), \D_n \in {\cal
D}_n {\rm \ and \ } n=k, k+1, \cdots \right\} \ \ , $$ is a
$\rho$--cover for $W_K(\p)$ with $\rho=\rho(k):= 3\,\p(2^k)$.

\medskip

\noindent{{\bf Step 3}: {\em The finale}. \ } Let $\rho=\rho(k):=
3\,\p(2^k)$.  Step 2 together with the  definition of
$s$--dimensional Hausdorff measure implies that $$ {\cal H}^s_\rho
(W_K(\p)) \ \leq \ \sum_{n=k}^{\infty} \ \ \  \sum_{\D_n \in {\cal
D}_n } \sum_{B_n(\p) \in {\cal C}(\D_n)} \!\! (3\, \p(2^n))^s \ \
. $$ Thus, in view of   (\ref{nod}) and (\ref{noc}),  it follows
that for $k$ sufficiently large  $$ {\cal H}^s_\rho (W_K(\p)) \
\ll \ \sum_{n=k}^{\infty} \ 2^{n \alpha \frac{d+1}{d} } \
\p(2^n)^{\alpha+s - \delta} \ . $$ This together with
(\ref{div2hd}) implies that $$ {\cal H}^s_\rho (W_K(\p)) \ \to \ 0
\hspace{7mm} {\rm as } \hspace{7mm} \rho \to 0 \ \ \ (k \to
\infty) \ \ , $$ and so ${\cal H}^s (W_K(\p)) = 0 $ as required.

\hfill $\spadesuit$
\newpage

\noindent {\bf Acknowledgements: \ }  SV would like to thank Barak
Weiss for  sending him the  pre-print \cite{KLW} from which this
work has materialized  and for numerous   e-mail conversations.
Also, thanks to Brent Everitt for  the    engaging conversations
regarding the $d$--simplex.  He would also like to thank Ayesha
and Iona for making 2003 such a fulfilling year -- may there be
many more. Much appreciated and thanks again  girls. Finally he
would like to thank the wonderful staff at `Bright Beginnings' who
have been absolutely fantastic with Ayesha and Iona -- in
particular their best friend Claire Kerner and Angela Foley for
taking them back.



\vspace{10mm}

\noindent Andy Pollington:  Department of Mathematics, Brigham
Young University,

 ~ \hspace{21mm}  Provo,   UT 84602,  USA.

\vspace{2mm}

 ~ \hspace{21mm} e-mail: : andy@math.byu.edu

\vspace{8mm}

\noindent Sanju L. Velani: Department of Mathematics, University
of York,

 ~ \hspace{19mm}  Heslington, York, YO10 5DD, England.

 \vspace{2mm}

 ~ \hspace{19mm} e-mail: slv3@york.ac.uk


\begin{thebibliography}{5}

\bibitem{falc}
K. Falconer : {\em Fractal Geometry: Mathematical Foundations and
Applications.} John Wiley \& Sons, (1990).



\bibitem{KLW}  D. Kleinbock, E. Lindenstrauss and B. Weiss :
On fractal measures and Diophantine approximation. {\em Selecta
Mathematica.}  To appear.


\bibitem{WT} B. Weiss : Almost no points on a Cantor set are very
well approximable. {\em Proc. R. Soc. London.} 457 (2001),
949--952.



\bibitem{W} B. Weiss : Dynamics on parameter spaces: submanifold and fractal
subset questions, in {\em Rigidity in Dynamics and Geometry}: M.
Burger and A. Iozzi (eds.). Springer, (2002).


\end{thebibliography}
\end{document}